\documentclass[conference]{IEEEtran}
\IEEEoverridecommandlockouts
\usepackage{upgreek}
\usepackage{cite}
\usepackage{amsmath,amssymb,amsfonts}
\usepackage{algorithmic}
\usepackage{graphicx}
\usepackage{textcomp}
\usepackage{xcolor}
\usepackage{comment}
\usepackage{mathtools}
\usepackage{amsmath, amssymb}
\def\BibTeX{{\rm B\kern-.05em{\sc i\kern-.025em b}\kern-.08em
    T\kern-.1667em\lower.7ex\hbox{E}\kern-.125emX}}
\begin{document}

\title{Simulation-Based Optimization for Policy Incentives and Planning of Hybrid Microgrids\\

}

\author{\IEEEauthorblockN{Nanrui Gong}
\IEEEauthorblockA{\textit{Department of Applied Mathematics and Statistics} \\
\textit{The Johns Hopkins University}\\
Baltimore, United States \\
ngong1@jh.edu}
\and
\IEEEauthorblockN{ James C. Spall}
\IEEEauthorblockA{\textit{Applied Physics Laboratory} \\
\textit{The Johns Hopkins University}\\
Applied Physics Laboratory\\
Laurel, Maryland 20723-6099 \\
james.spall@jhuapl.edu}

}

\maketitle

\begin{abstract}
Transitioning to renewable power generation is often difficult for remote or isolated communities, due to generation intermittency and high cost barriers. Our paper presents a simulation-based optimization approach for the design of policy incentives and planning of microgrids with renewable energy sources, targeting isolated communities.
We propose a novel framework that integrates stochastic simulation to account for weather uncertainty and system availability while optimizing microgrid configurations and policy incentives.
Utilizing the mixed-variable Simultaneous Perturbation Stochastic Approximation (MSPSA) algorithm, our method demonstrates a significant reduction in Net Present Cost (NPC) for microgrids, achieving a 68.1\% reduction in total costs in a case study conducted on Popova Island.
The results indicate the effectiveness of our approach in enhancing the economic viability of microgrids while promoting cleaner energy solutions.
Future research directions include refining uncertainty models and exploring applications in grid-connected microgrids.
\end{abstract}


\section{Introduction}

The Paris Agreement has defined a climate goal of no more than 1.5°C of global warming. To reach the climate goal, it is estimated that humans have to reach net zero by 2050 relative to a 2019 baseline \cite{UNNetZeroCoalition}. However, reducing emissions for geographically isolated power users remains a challenge. Geographic isolation prevents access to the main grid, keeping local users to burning fossil fuels for power generation. To assist isolated communities in transitioning towards cleaner sources of power generation, We propose a simulation-based optimization approach to jointly optimize policy incentives and planning of microgrids with renewable generation.    

Although distributed renewable generation have been proposed as an alternative power source for isolated communities, there are two main barriers to the adoption of renewable sources. The first barrier is high costs of adopting renewable sources. Isolated communities in developing countries are particularly sensitive to costs due to a shortage of investmentable funds. The second barrier is the uncertainties associated with intermittent power sources. The two most prominent forms of intermittent power generation, wind and solar, both suffer from unpredictable meteorological variations. The second barrier has been a major obstacle for the adoption of renewable sources in power grids world wide, as it leads to the risks of blackouts. 

To address the first barrier, effective microgrid sizing must be performed, accounting for weather uncertainties. To overcome the cost barrier, public policy incentives are often necessary to encourage CO2 emissions reduction. While previous studies have utilized stochastic optimization techniques for microgrid planning \cite{Barelli2020, He2023}, our method is unique in using simulation as a proxy for real systems. This enables us to consider the interactions between policy incentives and system design under uncertainty. By integrating these elements, our approach generates the optimal combination of policy incentives and component sizes that minimizes the Net Present Cost (NPC) of the intended microgrid.

Numerous previous studies have paid attention to microgrid sizing problem under  uncertanties, but don't consider the inovlvement of policy incentives.  Uwineza \cite{UWINEZA2021100607} uses the HOMER software to find lowest cost  configurations of microgrid for an isolated community, then uses Monte Carlo Simulation to model system uncertainties, but simulation and optimization are performed differently. Fioriti et al. \cite{9248897} uses MCS to model uncertainties of the microgrid and particle swarm optimization to minimize the NPC, but fails to consider system reliability. Wang et al. \cite{Wang2023} uses MCS simulation for both solar photovoltaics and wind power generation uncertainty as well as system reliability, but the MCS simulation and the snake optimization (SO) utilized in the paper are performed separately. 
A techno-economic analysis is carried out in \cite{Amini2022} that discusses the microgrid sizing problem under different policy incentives, but fails to take into account uncertainty of generation and system reliability.

The main contributions of our paper are follows:  
\begin{itemize}
\item Our paper proposes a novel minimum-cost model for designing policy incentives for renewable generation microgrids.
\item We incorporate a stochastic simulation model that account for both fluctuations in renewable energy generation and system component failures. The model is based on the Probabilistic Resource Adequacy Suite (PRAS) developed by National Renewable Energy Laboratory (NREL).
\item We optimize the noisy objective function using simulation-based optimization.
\end{itemize}

\section{Microgrid model}
\subsection{Solar PV model}
In energy systems analysis and modelling, power of solar photovoltaics are commonly modelled as a function of solar irrandiance and temperature. Therefore, we calculate the power of solar PVs following the convention of \cite{su11030683}:

\begin{equation}
P_{PV}=P_{STC}\cdot \frac{G_c}{G_{STC}} \cdot \left[1 + k_c \cdot (T_c - T_{STC})\right]\label{eq}
\end{equation}
 where $(P_{PV}$ (kW) is the actual power of the PV array. $P_{STC}$ (kW): Rated power under Standard Test Conditions (STC), $G_c$  (kW/m\(^2\)) is the solar irradiance striking the photovoltaic array. $k_c$ (\%/°C) is the temperature coefficient of power, $T_c$ (°C) is the photovoltaic temperature, and $T_{STC}$ is the ambient temperature under STC (25°C).

Irradiance is the uncertain component of the solar power model. From previous literature, solar irradiance of any given time point roughly follows a Gaussian distribution. For any given microgrid, we can then use a probabilistic TMY (Typical Meteorological Year) model to generate the profile of solar PVs \cite{Luo2022}. With knowledge of the variation of solar irradiance at a given site, we add independent Gaussian noises to the TMY, forming a stochastic profile of a year. Using this approach, we can best conserve the seasonal variations and daily variations of solar irradiance in the newly generated profile.

\subsection{Wind power model}
Using \cite{UWINEZA2021100607}, we model wind power generation as a function of wind speed:
\[
P_{WT} =
\begin{cases} 
0 & \text{if } v \leq v_{ci} \text{ or } v > v_{co}, \\
P_{\text{rated,WT}} \cdot \frac{v - v_{ci}}{v_{r} - v_{ci}} & \text{if } v_{ci} < v \leq v_{r}, \\
P_{\text{rated,WT}} & \text{if } v_{r} < v \leq v_{co}.
\end{cases}
\]
where $P_{WT} $ (W) is the power output of the wind turbine,
$v $(m/s) is the wind speed at the location of the wind turbine,
$v_{ci}$ (m/s) is the cut-in wind speed, below which the turbine does not generate power,
$ v_{co} $ (m/s) is the cut-out wind speed, above which the turbine shuts down to avoid damage, $ v_r$  (m/s) is the rated wind speed, at which the turbine generates its rated power.
$ P_{\text{rated, WT}}$ (W) is the rated power of the wind turbine, the maximum power output the turbine can generate at the rated wind speed. 

Fluctuations in wind speed is the main reason for causing volatility in wind power generation. To simulate flutuations in wind speed, we can generate wind speed data from a Weibull distribution \cite{SALGADODUARTE2020106237}. The parameters can be estimated based on given data. 

\subsection{Microturbine model}\label{SCM}

A microturbine generator is also included as a backup generator. We assume that the microturbine can serve all load that are not served by solar PV, wind turbines and the battery up to its maximum capacity (rated capacity). We also assume that the carbon emissions and fuel consumption are both proportional to the energy generated by the microturbine. 

\subsection{Battery storage model}
Following \cite{en15186706}, we define the battery storage model as follows. Let $\Delta t$ be the 1-hour simulation step, 
\begin{align}
S_{\text{}}(h) &= \eta_{\text{co}} S_{\text{}}(h-1) \nonumber \\
&\quad + \frac{P_{\text{ch}}(h) \Delta t \eta_{\text{ch}}}
{E_{\text{max}}^{\text{}}} H_{\text{ch}}(h) \nonumber \\
&\quad - \frac{P_{\text{dch}}(h) \Delta t}
{E_{\text{max}}^{\text{}} \eta_{\text{dch}}} H_{\text{dch}}(h)
\end{align} where $S_{\text{}}(h)$ is the state of charge (SOC) of the energy storage system (ESS) at the $h$-th moment,   $P_{\text{ch}}(h)$ and $P_{\text{dch}}(h)$ are the charging and discharging powers, respectively. $\eta_{\text{co}}$, $\eta_{\text{ch}}$ and $\eta_{\text{dch}}$ are  the carry-over, charging and discharging efficiencies, respectively. $E_{\text{max}}^{\text{}}$ is the max capacity of the battery, and $H_{\text{ch}}(h)$ and $H_{\text{dch}}(h)$ are the state functions of charging and discharging, respectively.

\subsection{Low-emissions policy incentives}
To encourage carbon emissions reductions in the microgrid and increase the use of renewable sources for power generation, a renewable energy penetration incentive and a carbon emissions reduction incentive are considered. The renewable energy penetration incentive is determined as a subsidy for the investment cost of the microgrid. The tuning of this policy incentive is considered as discrete scenarios in \cite{en15218285}, but here we formulate it as a continuous variable:
\[
I_{rp} =
\begin{cases} 
\sum \text{CAPEX}_i \cdot T_{rp} & \text{if } R_{rp} \geq T_{rp}, \\
0 & \text{if } R_{rp} < T_{rp}.
\end{cases}
\]
where $I_{rp}$ (USD) is the renewable energy penetration subsidy, $\text{CAPEX}_i$ (USD) is the investment cost of the $i$th component, $T_{rp}$ (\%) is the renewables penetration requirement threshold, and $R_{rp}$ (\%) is the actual renewables penetration of the current configuration. The renewable penetration subsidy only applies if the renewable penetration rate is at least the level of the rate determined by the policy. 

The carbon emission reduction incentive is also designed as a subsidy, as previously discussed in \cite{HE2023125454}. We model the scale of the subsidy as a portion of the carbon tax:
\[
I_{er} = 
\begin{cases} 
C_{tax} \cdot T_{er} & \text{if } T_{er} \geq R_{er}, \\
0 & \text{if } T_{er} < R_{er}.
\end{cases}
\]
where $I_{er}$ (\%) is the carbon emissions reduction subsidy for the microgrid, representing the financial incentive given for emissions reductions exceeding the required threshold. $C_{tax}$ (USD) is the total carbon tax cost paid by the microgrid in the baseline scenario where all load is served by microturbines (MT), $T_{er}$ (\%) is the 
 carbon emissions reduction requirement threshold that must be exceeded for the subsidy to apply, and $R_{er} $(\%) is the ctual emissions reduction rate achieved by the microgrid, expressed relative to the baseline case where all load is served by MT.

This subsidy only applies if the emissions reduction rate is at least the level of the rate determined by the policy. The emissions reduction rate is calculated as the ratio between the actual total carbon emissions and the hypothetical total carbon emissions if all load is served with the MT.
\subsection{Markov system reliability model}
Power system reliability is well studied using an availability model. Referencing \cite{SALGADODUARTE2020106237}, we use a two-state Markov process to simulate the availability and unavailibility of microgrid components.
\begin{align}
P_A = \frac{\mu}{\lambda + \mu}
\\
P_U = \frac{\lambda}{\lambda + \mu}
\end{align}
where $P_U$ is the stationary probability of component unavailability, and $P_A$ is the stationary probability of component availability. $\lambda$ is the failure rate and $\mu$ is the repair rate. The availability simulations are performed using the PRAS simulator \cite{stephen2021pras}. 

\subsection{Objective function}
Our objective function has the form:
\begin{align}
\text{NPC} 
&\coloneqq
\Bigg(\sum \text{CAPEX}_i + \sum \text{OPEX}_i \\
&\quad + C_{\text{tax}} + C_{\text{fuel}} + C_{\text{voll}} \\
&\quad - I_{\text{rp}} - I_{\text{er}} \Bigg) \cdot \text{CRF}
\end{align}
subject to the constraints:
\begin{align}
P_i^{\text{min}} \leq P_i \leq P_i^{\text{max}}
\\
T_i^{\text{min}} \leq T_i \leq T_i^{\text{max}}
\\
\text{HLL}\leq h_{max}
\end{align}where NPC (USD) is the Net Present Cost, $\sum \text{OPEX}_i$ (USD) is the total operational and maintenance costs of the microgrid, $C_{\text{voll}}$ (USD) is the cost of lost load, and CRF is the cost recovery factor. $P_i^{\text{min}}$  (kWh) and $P_i^{\text{max}}$  (kWh) are the lower and upper bounds of the component rated capacity, respectively. $T_i^{\text{min}}$ and $T_i^{\text{max}}$ are the lower and upper bounds of the policy incentives requirement threshold, respectively. \text{HLL} (hr) is the total hours of lost load across the year, and  $h_{max}$ (hr) is the total hours of lost load allowed.

\section{Constrained MSPSA Algorithm}
We select the mixed simultaneous perturbation stochastic approximation (MSPSA) \cite{8430974,Wang2025}, given that the problem has both discrete and continuous variables. The mixed variable version is developed from the general simultaneous perturbation stochastic
approximatio (SPSA) and discrete SPSA methods \cite{Spall1992, Spall1994,Wang2013}.  
SPSA has seen applications in other resource allocation problems in \cite{Chin1999}.

\begin{enumerate}
    \item   
    Set counter index $k=0$. Select initial guess vector $\boldsymbol{\hat{\uptheta}}_{0}$ and nonnegative coefficients $a$, $c$, $A$, $\alpha$, and $\gamma$ in the MSPSA gain sequences
    $
    a_k = a/(k + 1 + A)^\alpha$ and $c_k = c/(k + 1)^\gamma.$
    
    \item  
    Generate a $p$-dimensional random perturbation vector $\Delta_k$ using a Monte Carlo algorithm, where each of the $p$ components of $\Delta_k$ is independently generated from a mean-zero probability distribution satisfying certain regularity conditions (discussed later).  
    An effective (but not mandatory) choice is to use a Bernoulli $\pm 1$ distribution with probability $1/2$ for each outcome.

    \item  
    Obtain two measurements of the loss function based on the simultaneous perturbation around the current estimate $\boldsymbol{\hat{\uptheta}}_{k}$:  
    \[
    y_k^{(+)} \quad \text{and} \quad y_k^{(-)},
    \]
    using $\mathbf{C}_k$ and $\Delta_k$ from Steps 1 and 2.

    \item   
    Generate the simultaneous perturbation gradient approximation $\hat{\mathbf{g}}_k(\boldsymbol{\hat{\uptheta}}_{k})$ using 
    \[\hat{g}_k(\boldsymbol{\hat{\uptheta}}_{k}) = \frac{ y_k^{(+)} -  y_k^{(-)}}{2 C_k \odot \Delta_k}
    \]
    where \( (C_k \odot \Delta_k)^{-1} \) is the vector of inverses of the components of \( C_k \odot \Delta_k \) and $\odot$ is the matrix Hadamard product.
    \item  
    Use the standard SPSA form to update $\boldsymbol{\hat{\uptheta}}_{k}$ to a new value $\boldsymbol{\hat{\uptheta}}_{k+1}$.
\[
\boldsymbol{\hat{\uptheta}}_{k+1}= \boldsymbol{\hat{\uptheta}}_{k} - a_k \hat{\mathbf{g}}_k( \boldsymbol{\hat{\uptheta}}_{k})
\]
    \item   
    Return to Step 1 with $k + 1$ replacing $k$. Terminate the algorithm if there is little change in several successive iterates or if the maximum allowable number of iterations has been reached.  
    Return the final estimate $\boldsymbol{\hat{\uptheta}}_{N}$ where the rounding operator $[\cdot]$ acts on all discrete variables.
\end{enumerate}

To implement constraints on the SPSA algorithm, we use a simple projection technique:

\[
\boldsymbol{\hat{\uptheta}}_{k+1}= \mathcal{P}\left(\boldsymbol{\hat{\uptheta}}_{k} - a_k \hat{\mathbf{g}}_k( \boldsymbol{\hat{\uptheta}}_{k})\right)
\]where $\mathcal{P}(\boldsymbol{\hat{\uptheta}}_{k})$ is the nearest point in the feasible region to $\boldsymbol{\hat{\uptheta}}_{k}$ in terms of Euclidean distance.

Constraint (10) is both stochastic and unknown, and thus cannot be easily implemented by projection. To enforce the constraint, we add a quadratic penalty term to the objective function \cite{1271742}. The objective function is thus modified as: 
\begin{equation}
L(\uptheta) = NPC + r\left[\max\{0, q(\uptheta)\}\right]^2 
\end{equation}where $q(\uptheta)=\text{HLL}-h_{max}$. $r$ is a positive real number normally referred to as the penalty parameter. $r$ should be carefully tuned through trial and error to ensure convergence.

\section{Case Study}

In this section, we apply our algorithm to the Popova Island, located near Vladivostok in the Pacific. We derive our data from \cite{UWINEZA2021100607}. Since we only have meterological data of a single year 2013, we use the single year data in place of the TMY for solar irradiance simulation. The average standard deviation is $\sigma_{PV}=72.4$. We then fit a Weibull distribution to aggregated wind speed across the year, and attain shape and scale parameters, $k_{WT}=2.0$, $\lambda_{WT}=6.8$, respectively. The distribution parameters are used to generate stochastic solar and wind power profile across the year. 
Parameters for the MSPSA algorithm are chosen as $\alpha=0.602, A=500, a=0.25, c=0.7, \gamma=0.101$. Choices of the parameters are based on repeated tests on case study data to ensure stability and convergence.   
We set the initial vector as \[\boldsymbol{\hat{\uptheta}}_{0} = [5000,5000, 5000, 5000, 0.0, 0.0]^T
\]
The reason for choosing the above $\boldsymbol{\hat{\uptheta}}_{0}$ is so that we can start from a near-midpoint of the component capacity bounds with no support from policy incentive. 
To assess the performance of the MSPSA optimization algorithm, we compare it with a mixed variable particle swarm optimization approach. PSO has been widely used for microgrid planning problems \cite{RECALDE2020106540}. Parameters of PSO include cognitive coefficient $c_1=2.3$, social coefficient $c_2=2.3$, inertia weight $w=1$, initial velocity $v_0 \sim U(-1,1)$, and population size $m=20$. We set the initial position of PSO as $\boldsymbol{\hat{\uptheta}}_{0}$ to ensure a fair comparison. Detailed formulation of the algorithm is omitted here.

\section{Results}
 
\begin{figure}[htbp]
\centerline{\includegraphics[scale=0.40]{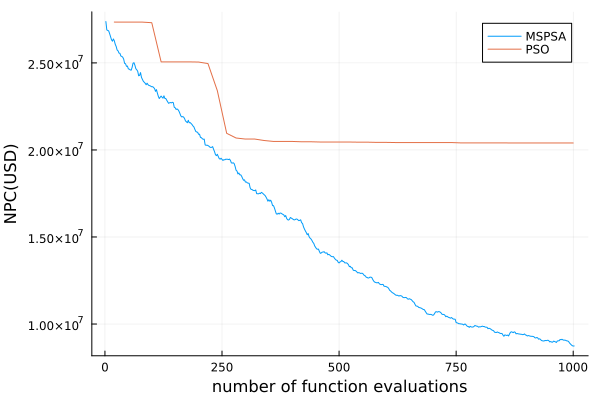}}
\caption{Comparison of MSPSA and PSO performance over 10 replicates}
\label{fig}
\end{figure}

 After 500 iterations, MSPSA achieves a 68.1\% reduction in total cost, compared to the initial state (Fig.1). In comparison, the PSO algorithm achieves a reduction 25.4\% in objective function value with the same amount of function evaluations. Both MSPSA and PSO convergence curves are attained after averaging over 10 replicates. Therefore, MSPSA has clear advantage in convergence rate in our case study. Function value fluctuations during the optimization steps are not pronounced in our case likely due to a lower level of randomness in the simulation. 
 \begin{figure}[htbp]
\centerline{\includegraphics[scale=0.40]{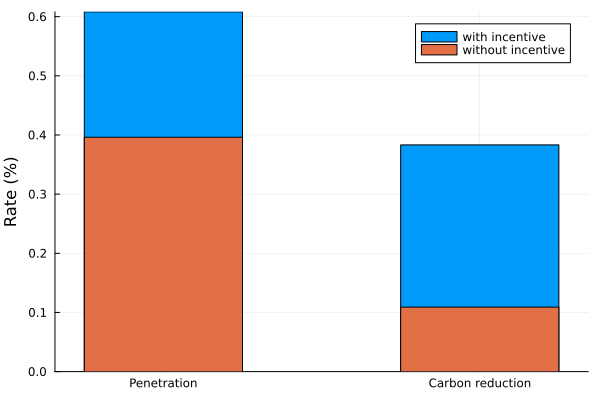}}
\caption{Comparison of planning outcome with and without policy incentives}
\label{fig}
\end{figure}

\begin{table}[htbp]
\caption{Comparison of microgrid configurations}
\centering
\begin{tabular}{|c|c|c|c|c|}
\hline
\textbf{Scenarios} & \multicolumn{4}{|c|}{\textbf{Solutions}} \\
\cline{2-5} 
\textbf{}  & 
\textbf{\textit{PV} (kW)} & \textbf{\textit{WT} (kW)} & \textbf{\textit{BSS} (kWh)} & \textbf{\textit{MT} (kW)}  \\
\hline

Without incentives & 1380	
  & 103	 & 3064	   & 2962 \\
\hline
With incentives & 1686	
  & 1077 &1783&1756 \\

\hline
\end{tabular}
\label{tab1}
\end{table}

Based on the solution of MSPSA at the 500th iteration during one replication, we obtained a current renewable energy penetration requirement of 22.1\% and a carbon reduction requirement of 6.9\%. To evaluate the effectiveness of the policy incentives, we set the policy incentive variables to zero and optimized again. The planning results for the scenarios with and without incentives are described in Table 1. When evaluating both scenarios using a Monte Carlo simulation with 100 samples, the average renewable energy penetration rate and carbon emission reduction rate were significantly higher when policy incentives were applied. The average NPC was also lower when policy incentives were applied, but this outcome was expected since the optimization algorithm was free to explore the case where the policy incentive variables were set to zero. The interpretation of the current solution is that applying a combination of policy incentives shifts the most economical microgrid configuration toward higher renewable generation and reduced reliance on fossil fuels, leading to greater carbon reduction. However, since we are not certain whether the solution at the 500th iteration is globally optimal, this interpretation should not be generalized.  


\section{Conclusion}
This paper presents a novel framework for jointly optimizing policy incentives and planning  microgrid. We proposed a simulation-based optimization framework using the mixed-variable MSPSA algorithm to address the dual challenges of renewable energy uncertainties and high upfront costs. Our approach incorporates stochastic simulations to account for renewable generation variability, system reliability, and component failures. Furthermore, the framework optimizes the combinations of policy incentives, such as renewable penetration subsidies and carbon emissions reduction subsidy, to enhance economic viability and reduce emissions.

In the case study of Popova Island, our MSPSA-based method demonstrated superior performance in cost reduction and convergence efficiency compared to the conventional PSO algorithm. MSPSA achieved an average 68.1\% reduction in total cost compared to a hypothetical initial case, highlighting its potential as an effective tool for addressing the challenges of hybrid microgrid planning under uncertainty. Further, comparing renewable energy penetration rate and carbon emission reduction rate with and without policy incentives suggest the usefulness of our framework for tuning low-carbon policy incentives.

Our study utilizes a fairly simple model for simulating power system uncertainties. Future work should combine the MSPSA approach with more realistic uncertainty simulation models, including the correlation between weather and power usage variations. In addition, another future direction involves applying our algorithm to grid-connected microgrids. By bridging simulation, optimization, and policy design, this study contributes to sustainable energy solutions for isolated communities and supports global efforts toward carbon neutrality.

\bibliographystyle{IEEEtran}
\bibliography{references}

\vspace{12pt}
\color{red}

\end{document}